\documentclass[10pt,a4paper]{amsart}
\usepackage{amssymb, amsmath, amsthm}

\theoremstyle{plain}
\newtheorem{thm}{Theorem}

\newtheorem{lem}{Lemma}

\theoremstyle{definition}

\theoremstyle{remark}
\newtheorem{rem}[subsubsection]{Remark}
\newtheorem{exam}[subsubsection]{Example}
\newtheorem*{prob}{Problem}
\newtheorem*{claim}{Claim}

\DeclareMathOperator{\dist}{dist}
\DeclareMathOperator{\Log}{Log}

\newcommand{\imax}{\underline{m}}

\title[Stability of Analytic Germs]{On the Stability of Analytic Germs under Ultradifferentiable Perturbations}

\author{Vincent Thilliez}

\address{Math\'ematiques - B\^atiment M2\\
Universit\'e des Sciences et Technologies de Lille\\
F-59655 Villeneuve d'Ascq Cedex, France}

\email{thilliez@math.univ-lille1.fr}


\begin{document}

\begin{abstract} Let $ f$ be a real-analytic function germ whose critical locus contains a given real-analytic set $ X $, and let $ Y $ be a germ of closed subset of $ \mathbb{R}^n $ at the origin. We study the stability of $ f $ under perturbations $ u $ that are flat on $ Y $ and that belong to a given Denjoy-Carleman non-quasianalytic class. We obtain a condition ensuring that $ f+u=f\circ\Phi $ where $ \Phi $ is a germ of diffeomorphism whose components belong to a (generally larger) Denjoy-Carleman class. Roughly speaking, this condition involves a \L ojasiewicz-type separation property between $ Y $ and the complex zeros of a certain ideal associated with $ f $ and $ X $. The relationship between the Denjoy-Carleman classes of $ u$ and $ \Phi $ is controlled precisely by the inequality. This result extends, and simplifies, former work of the author on germs with isolated critical points.   
\end{abstract}

\maketitle

\section*{Introduction}

The stability of a smooth function germ under flat perturbations can be expressed by the notion of infinite determinacy, whose most basic case goes as follows. Denote by $ \mathcal{E}_n $ the ring of $ C^\infty $ function germs
at the origin in $ \mathbb{R}^n $, and by $ \imax^\infty $ the ideal of flat germs, that is, germs vanishing at infinite order at the origin. An
element $ f $ of $ \mathcal{E}_n $ is \emph{infinitely determined} if, for
any element $ u $ of $ \imax^\infty $, there exists a germ $ \Phi $ of
$ C^\infty $-diffeomorphism at the origin such that $ f+u=f\circ\Phi $. In abbreviated form, this can be written
\begin{equation}\label{deter00}
f+\imax^\infty\subset f\circ\mathcal{R}, 
\end{equation}
where $ \mathcal{R} $ denotes the group of
germs of $ C^\infty $-diffeomorphisms at the origin. Property \eqref{deter00} can be characterized by means of the ideal
$ \langle\nabla f\rangle\mathcal{E}_n $ generated in $ \mathcal{E}_n $ by the partial
derivatives $ \partial f/\partial x_j $, $ j=1,\ldots,n $. Indeed, $ f  $ is infinitely determined 
if and only if
\begin{equation}\label{jac00}
\imax^\infty\subset \langle\nabla f\rangle\mathcal{E}_n,
\end{equation}
see \cite{Nguyen, Wilson} or part II of \cite{Wall}. 

The case of real analytic germs is particularly important since, in this case, condition \eqref{jac00} holds if and only if $ f $ has at most an isolated real critical point at $ 0 $. For such germs, we studied in \cite{Thilliez1} a quantitative version of infinite determinacy involving perturbations and diffeomorphisms chosen in suitable subclasses of $\mathcal{E}_n $ and $ \mathcal{R} $, namely non-quasianalytic ultradifferentiable classes. The precise definition of these
classes will be recalled in subsection \ref{classes}.
For the moment, let us only mention the typical Gevrey example: a smooth germ belongs to the Gevrey 
ring $ \mathcal{C}_{n,\alpha} $, where $ \alpha $ is a positive real number, if, in a suitable neighborhood of $ 0 $, its derivatives at any order $ j $ are bounded by
$ C^{j+1}j!^{1+\alpha} $ for some positive constant $ C $. Put $
\imax^\infty_\alpha=\imax^\infty\cap\mathcal{C}_{n,\alpha} $ and denote by
by $ \mathcal{R}_\alpha $ the set of those elements of $ \mathcal{R} $ whose components belong to 
$ \mathcal{C}_{n,\alpha} $. Now let $ f $ be a real-analytic function germ with an
isolated real critical point at $ 0 $. 
Theorem 4.1 of
\cite{Thilliez1} asserts that
\begin{equation}\label{isgev}
f+\imax^\infty_\alpha\subset f\circ\mathcal{R}_\beta,
\end{equation}
where the ratio $ \beta/\alpha $ equals a suitable \L ojasiewicz exponent related to the set of
{\it complex} critical points of $ f $. In general, such a loss of regularity is unavoidable, and the value of $ \beta/\alpha $ given by \cite{Thilliez1} is sharp. For instance, when $ n=2 $ and
$ f(x)=(x_1^2+x_2^4)^2 $, it is possible to take $ \beta=2\alpha $ but not any smaller $ \beta $. In particular, if $ u $ denotes the element of
$ \imax^\infty_\alpha $ given by $ u(x)=x_1^2(x_1^2+x_2^4)\exp(-\vert x_2\vert^{-1/\alpha}) $, it is shown in \cite{Thilliez1} that $ f+u $ does not belong to $ f\circ\mathcal{R}_\beta $ for $ \beta<2\alpha $. The factor
$ 2 $ is precisely the \L ojasiewicz exponent for the
regular separation between the real plane $ \mathbb{R}^2 $ and the complex critical set 
$ \{x\in\mathbb{C}^2 : x_1^2+x_2^4=0\} $.

The study of infinite determinacy for non-isolated singularities began more recently in the $ C^\infty $ setting \cite{Grandjean2, Sun-Wilson, Thilliez3}. A common framework for these papers can be described as follows. Consider 
a real-analytic map germ
$ \psi : (\mathbb{R}^n,0)\rightarrow (\mathbb{R}^p,0) $. Denote by $ \langle\psi\rangle $ the ideal generated by the components of $ \psi $ in the ring $ \mathcal{O}_n $ of real-analytic function germs, and by $ \int\langle\psi\rangle $ its primitive ideal, that is, the ideal of elements of $ \langle\psi\rangle $ whose first order derivatives also belong to $ \langle \psi\rangle $. One studies elements $ f$ of $ \int\langle\psi\rangle $ (see subsection \ref{fKf} for the motivation). Using a classical homotopy technique, one can then show (see remark \ref{smoothcase} below) that the condition 
\begin{equation}\label{jacX0}
\langle\psi\rangle\imax^\infty\subset\langle\nabla f\rangle\mathcal{E}_n 
\end{equation}
implies the determinacy property
\begin{equation}\label{deterX0}
f+\int\langle\psi\rangle\imax^\infty\subset f\circ\mathcal{R}.
\end{equation}
Obtaining a converse implication is much more difficult, and requires the replacement of \eqref{deterX0} by a more precise property taking into account the preservation of $ \psi^{-1}(\{0\}) $ by diffeomorphisms, as well as additional properties of $ \psi $, see \cite{Thilliez3}. We shall not consider this question here (see, however, remark \ref{presX} below).

As in the case of isolated singularities, the implication \eqref{jacX0}$\Longrightarrow$\eqref{deterX0} is particularly interesting when $ f $ is analytic: indeed, \eqref{jacX0} is then equivalent to a condition of real isolated zero for a suitable Fitting ideal $ \mathcal{K}_f $ associated with $ \psi $ and $ f $. Accordingly, the purpose of the present paper is twofold. First, we shall extend the aforementioned theorem 4.1 of \cite{Thilliez1} to the case of real-analytic germs with non-isolated singularities. Second, we will study an 
aspect of determinacy which is not described in \cite{Thilliez1}, namely, how regularity estimates can be influenced by the flatness of perturbations along certain sets. In the example of loss of regularity given above, the perturbation $ u $ is flat on the $ x_1 $-axis. If we consider, instead, any perturbation $ v $ that is flat on the $ x_2 $-axis, it turns out\footnote{See subsection \ref{examp} for more details.} that $ f+v $ belongs to $ f\circ\mathcal{R}_\alpha $: in other words, for such a $ v $, there is no loss of regularity. Of course, this phenomenon is hidden in the $ C^\infty $ case, for which the idea
of assuming flatness on various subsets of $ \mathbb{R}^n $ has, nevertheless, been used in \cite{Kushner-Terra Leme, Thilliez3}. 

Thus, in order to extend \eqref{isgev} with respect to both aspects,  
we shall study here the condition
\begin{equation}\label{nonisgev}
f+\int\langle\psi\rangle\imax^\infty_{Y,\alpha}\subset f\circ\mathcal{R}_\beta,
\end{equation}
where $ Y $ is a germ of closed subset of $\mathbb{R}^n $ and $ \imax^\infty_{Y,\alpha} $ denotes the ideal of germs of $\mathcal{C}_{n,\alpha} $ that are flat on $ Y $. The goal is then to obtain a sharp sufficient condition for \eqref{nonisgev} in terms of complex zeros of elements of $ \mathcal{K}_f $. This will lead to theorem \ref{main} below.

Here is an example in $ \mathbb{R}^3 $. Consider $ f(x)=(x_1^2+x_2^2+x_3^4)^2(x_1^2+x_2^2) $. Put $ Y=\{x\in\mathbb{R}^3 : \vert x_1\vert=\vert x_2\vert=\vert x_3\vert^\mu\} $, where $ \mu $ is a given positive real number. We shall obtain the determinacy property $ f+\langle x_1^2,x_1x_2,x_2^2\rangle\imax^\infty_{Y,\alpha}\subset f\circ\mathcal{R}_\beta $ with $ \beta=2\alpha $ if $ 0<\mu\leq 1 $, resp. $ \beta=2\alpha/\mu $ if $ 1<\mu\leq 2 $ and $ \beta=\alpha $ if $ \mu>2 $.

As usual with determinacy results, a Mather-type homotopy argument reduces the
proof of theorem \ref{main} to the construction of suitable vector
fields. This construction will be achieved by means of somewhat more algebraic arguments than in \cite{Thilliez1}. Incidentally, it yields an alternative, simpler, proof of theorem 4.1 of \cite{Thilliez1}, as explained in subsection \ref{isolsi}.

It should also be emphasized that the results of the paper are not limited to the Gevrey scale, although this is the only one mentioned in the introduction for sake of simplicity. We have in fact a greater degree of freedom, which allows to sharpen the estimates by means of some tools from \cite{Thilliez0}: for instance, in the preceding example, if $ Y=\{x\in\mathbb{R}^3 : \vert x_1\vert=\vert x_2\vert=\vert x_3\vert^\mu\Log(1+\frac{1}{\vert x_3\vert})^{\nu} \} $ with $ 1<\mu< 2 $ and $ \nu>0 $, then the Gevrey estimate $ j!^{1+\beta} $ for the elements of $ \mathcal{R} $ involved in the determinacy property becomes
$ j!^{1+2\alpha/\mu}(\Log j)^{2\nu j/\mu} $. 

\section{Preliminaries and a Division Result}\label{basic}

\subsection{Notation}\label{nota} Throughout the paper, we use the following notation. 

Let $ \mathcal{A} $ and $ \mathcal{B} $ be two subrings of $ \mathcal{E}_n $ with $ \mathcal{A}\subset\mathcal{B} $. If $ \phi=(\phi_1,\ldots,\phi_q) $ is a $ \mathbb{R}^q $-valued map germ with components in $ \mathcal{A} $, we denote by $ \langle\phi\rangle \mathcal{B} $ the ideal generated by $ \phi_1,\ldots,\phi_q $ in $ \mathcal{B} $. When $ \mathcal{A}=\mathcal{B}=\mathcal{O}_n $, we omit it in the notation and write simply $ \langle\phi\rangle $. If $ \mathcal{I}$ is an ideal of $ \mathcal{A} $, we denote by $ \int\mathcal{I} $ its primitive ideal, that is, the ideal of elements of $ \mathcal{I} $ whose first order derivatives also belong to $ \mathcal{I}$. For any $ \mathcal{B} $-module $ \mathcal{M} $, 
we denote by $ \mathcal{IM} $ the submodule generated by $ \mathcal{I} $ in $ \mathcal{M} $.

For any element $ f $ of $ \mathcal{E}_n $ and any subset $ \mathcal{R}' $ of $ \mathcal{R} $, the set
$ \{f\circ \Phi :\Phi\in \mathcal{R}'\} $ is denoted by $ f\circ\mathcal{R}' $.

For any multi-index $ J=(j_1,\ldots,j_n) $ in $ \mathbb{N}^n $, the length $ j_1+\ldots+j_n $ of $ J $ is denoted by the corresponding lower case letter $ j $. We
put $ D^J =\partial^j/\partial x_1^{j_1}\cdots\partial x_n^{j_n} $. The gradient of a smooth function with respect to $ x_1,\ldots,x_n $ is denoted by $ \nabla_x $, or simply by $ \nabla $ if there is no variable. 

Properties that hold on a given subset $ V $ of
$ \mathbb{R}^n $ containing the origin are always
understood in the sense of germs, that is, on some representative of $ V $.

\subsection{Moderate growth sequences and admissible functions}\label{sequ}
We recall here the basic properties of sequences that will be needed to deal with  ultradifferentiable classes. All the proofs and details can be found in \cite{Thilliez0} and the references therein. 

In what follows, we say that a sequence $ M=(M_j)_{j\geq 0} $ of real numbers is \emph{tame} if it satisfies the following conditions:
\begin{equation}\label{seq1}
\textrm{normalization: the sequence } M \textrm{ is increasing, with } M_0=1,
\end{equation}
\begin{equation}\label{seq2}
\textrm{logarithmic convexity: the sequence } (M_{j+1}/M_j)_{j\geq 0} \textrm{ is increasing,}
\end{equation}
\begin{equation}\label{seq3}
\textrm{moderate growth: there exists a constant } A>0 \textrm{ such that } 
\end{equation}
\begin{equation*}
M_{j+k}\leq A^{j+k}M_jM_k \textrm{ for any } (j,k)\in\mathbb{N}^2.
\end{equation*}

For any real $ t\geq 0 $, put now $ h_M(t)=\inf_{j\geq 0}t^jM_j $. The function $ h_M $ then determines $ M $ since the logarithmic convexity assumption implies $ M_j=\sup_{t>0}t^{-j}h_M(t) $.
Note that \eqref{seq1} and \eqref{seq2} imply
\begin{equation}\label{seq4} 
M_jM_k\leq M_{j+k} \textrm{ for any } (j,k)\in\mathbb{N}^2,
\end{equation} 
so that \eqref{seq3} amounts to a similar estimate in the opposite direction. Another consequence is the existence of a constant $ B>0 $ such that
\begin{equation}\label{hMsq}
h_M(t)\leq \big(h_M(Bt)\big)^2\textrm{ for any } t\geq 0.
\end{equation}

Now, let $ \theta $ be a continuous increasing function on $ [0,\varepsilon[ $ for some $ \varepsilon>0 $, with $ \theta(0)=0 $. We shall say that $ \theta $ is \emph{admissible} if $ t\mapsto \theta(t)/t $ is increasing on $ ]0,\varepsilon[ $ and if there exists $ s\geq 1 $ such that $ t\mapsto \theta(t)/t^s $ is decreasing on $ ]0,\varepsilon[ $.
The following result appears in \cite{Thilliez0}.

\begin{lem} Let $ M $ be a tame sequence and let $ \theta $ be an admissible function. Then there exists a tame sequence $ M^{(\theta)} $ for which one can find positive constants $ c $ and $ c' $ such that, for any sufficiently small $t>0$, 
\begin{equation*}
h_M(ct)\leq h_{M^{(\theta)}}(\theta(t))\leq h_M(c't). 
\end{equation*} 
This sequence is unique modulo the equivalence relation $ \sim $ defined as follows: $ M'\sim M'' $ if there exists a constant $ C\geq 1 $ such that $ C^{-j}M''_j\leq M'_j \leq C^jM''_j $ for any $ j $.
\end{lem} 

\begin{rem}\label{grtMtheta} 
It is not difficult to show that there is a constant $ C\geq 1 $ such that $ C^{-j}M_j\leq M^{(\theta)}_j\leq C^j(M_j)^s $ for any $ j $.
\end{rem}

\begin{exam}\label{expow} The most basic example occurs when $ \theta(t)=ct^s $ for given real numbers $ c>0 $ and $ s\geq 1 $. Indeed, we then have $ M^{(\theta)}_j=(M_j)^s $.
\end{exam} 

\begin{exam}\label{expowlog} Let $ \alpha $, $ \beta $, $ \mu $, $ \nu $ be real numbers with $ \alpha>0 $, $ \beta\geq 0 $, $ \mu>1 $ and $ \nu\geq 0 $. Put $ M_j=j!^\alpha\big(\Log(e+j)\big)^{\beta j} $. Then $ M $ is tame and we have $ H(at)\leq h_M(t)\leq H(bt) $ for $H(t)= \exp\big(-t^{-1/\alpha}\big(\Log(1/t)\big)^{-\beta/\alpha}\big) $ and for suitable positive constants $ a $ and $ b $. The function $ \theta(t)=t^\mu\big(\Log(1+\frac{1}{t})\big)^{-\nu} $ is admissible and we can take 
$ M^{(\theta)}_j=j!^{\alpha\mu}\big(\Log(e+j)\big)^{(\beta\mu+\nu)j} $.
\end{exam}

\subsection{Ultradifferentiable classes of germs}\label{classes}
 
A tame sequence $ M $ is said to be \emph{non-quasianalytic} if
\begin{equation}\label{nonqa}
\sum_{j\geq 0}\frac{M_j}{(j+1)M_{j+1}}<\infty. 
\end{equation}
Being given such a sequence, the Denjoy-Carleman class of germs $ \mathcal{C}_{n,M} $ is defined as the set of those elements $ u $ of
$ \mathcal{E}_n $ for which one can find a neighborhood $ U $ of $ 0 $ in $ \mathbb{R}^n $ and a positive constant $ C $, depending on $ u $, such that the estimate
\begin{equation*}
\vert D^Ju(x)\vert \leq  C^{j+1}j!M_j 
\end{equation*}
holds for every point $ x $ in $ U $ and every multi-index $ J $. The sequence
$ M $ measures, in some sense, the defect of analyticity of the elements of $ \mathcal{C}_{n,M} $. 
The assumptions on $ M $ ensure that $ \mathcal{C}_{n,M} $  
is a local ring, stable by composition and derivation, whose maximal ideal $ \imax_M $ is generated by the coordinate functions. The non-quasianalyticity condition \eqref{nonqa} ensures that the ideal $ \imax_M^\infty $ of flat germs
in $ \mathcal{C}_{n,M} $ is not reduced to $ \{0\} $. More generally, if $ Y $ is a germ of closed subset of $ \mathbb{R}^n $ at the origin, we denote by $ \imax^\infty_{Y,M} $ the ideal of germs in $ \mathcal{C}_{n,M} $ that are flat on $ Y $ (in particular, $ \imax_M^\infty=
\imax^\infty_{\{0\},M} $).

The following notation will be used in the proof of theorem \ref{main}:
for any real number $ a$, we denote by $ \mathcal{C}_{n,M}(a) $ the ring of germs of functions $ \widetilde{u} $ at $ (0,a) $ in $ \mathbb{R}^{n+1} $ such that $ (x,t)\mapsto\widetilde{u}(x,t-a) $ belongs to $ \mathcal{C}_{n+1,M} $ (in particular, we have $ \mathcal{C}_{n,M}(0)=\mathcal{C}_{n+1,M} $). The corresponding maximal ideal is denoted by $ \imax_M(a) $. Being given a germ $ Y $ of closed subset at the origin of $ \mathbb{R}^n $, we denote by $ \imax^\infty_{Y,M}(a) $ the ideal of those elements of $ \mathcal{C}_{n,M}(a) $ that are flat on the germ of $ Y\times \mathbb{R} $ at $ (0,a) $. 

\begin{rem}\label{injec}
There is a natural injection $ \mathcal{C}_{n,M}\rightarrow \mathcal{C}_{n,M}(a) $; indeed, any element $ u $ of $ \mathcal{C}_{n,M} $ can be identified with the element $ \widetilde{u} $ of $ \mathcal{C}_{n,M}(a) $ defined by $ \widetilde{u}(x,t)=u(x) $. This identification maps $ \imax^\infty_{Y,M} $ into $ \imax^\infty_{Y,M}(a) $.
\end{rem}

\subsection{Fitting ideals and a division result}\label{singuls}
Consider two real analytic map germs $ \psi:(\mathbb{R}^n,0)\rightarrow(\mathbb{R}^p,0) $ and $ \varphi:(\mathbb{R}^n,0)\rightarrow(\mathbb{R}^q,0) $ such that $ \langle\varphi\rangle\subset\langle\psi\rangle $. Put  
$ \mathcal{M}=\sigma^{-1}(\langle\varphi\rangle) $, where $ \sigma : \mathcal{O}_n^p\longrightarrow \mathcal{O}_n $ is the map given by 
$ \sigma(f_1,\ldots,f_p)=\sum_{i=1}^pf_i\psi_i $. Clearly, $ \mathcal{M} $ is the submodule of $ \mathcal{O}_n^p $ generated
by $ h^1,\ldots,h^q,k^1,\ldots,k^r $, where each $ h^j $ is an element of $ \mathcal{O}_n^p $ such that 
$ \varphi_j=\sigma(h^j) $ and $k^1,\ldots,k^r $ is a system of generators for the 
module of relations $ \ker\sigma $.  One can thus write 
$ \mathcal{M}=\lambda(\mathcal{O}_n^{q+r}) $, where 
$ \lambda : \mathcal{O}_n^{q+r}
\longrightarrow \mathcal{O}_n^p $ is the morphism of free modules given by
$\lambda(\xi_1,\ldots,\xi_{q+r})= \sum_{j=1}^q\xi_jh^j+\sum_{j=1}^r\xi_{q+j}k^j$.
Denote by $ \mathcal{K} $ the ideal generated in $ \mathcal{O}_n $ by the maximal minors of the
matrix of $ \lambda $ in the canonical bases of $ \mathcal{O}_n^{q+r}
$ and $ \mathcal{O}_n^p $. 
Following section 20.2 of \cite{Eisenbud}, we see that $ \mathcal{K} $ is precisely 
the Fitting ideal $ \mathsf{Fitt}_0\big(\mathcal{O}_n^p/\mathcal{M}\big) $. In particular, it depends only on $ \psi $ and $\varphi $, and we have the inclusion
$ \mathcal{K}\mathcal{O}_n^p\subset \mathcal{M} $.
Hence, taking the image by $ \sigma $, we get
\begin{equation}\label{dol}
\mathcal{K}\langle\psi\rangle\subset \langle\varphi\rangle. 
\end{equation}

The following result, which can be viewed as a division theorem, extends theorem 2.3 of \cite{Thilliez2} (dealing with the particular case $ \psi(x)=x $ and $ Y=\{0\} $).

\begin{thm}\label{divis}
Let $ \psi:(\mathbb{R}^n,0)\rightarrow(\mathbb{R}^p,0) $ and $ \varphi:(\mathbb{R}^n,0)\rightarrow(\mathbb{R}^q,0) $ be two real-analytic map germs such that 
\begin{equation*}
\langle\varphi\rangle\subset\langle\psi\rangle. 
\end{equation*}
Let $ Y $ be a germ of closed set at the origin of $ \mathbb{R}^n $ and let $ \theta $ be an admissible function. Assume that the ideal $ \mathcal{K} $ contains a germ $ g $ whose complex zero set $ Z_g $ satisfies the generalized \L ojasiewicz inequality
\begin{equation}\label{lojagen}
\dist(x,Z_g)\geq \theta\big(\dist(x,Y)\big)\text{ for any real point }x.
\end{equation}
Then, for any non-quasianalytic tame sequence $ M $, we have 
\begin{equation*}
\langle\psi\rangle\imax^\infty_{Y,M}\subset\langle\varphi\rangle\imax^\infty_{Y,M^{(\theta)}}.
\end{equation*}
\end{thm}

\begin{proof} Since a part of the proof goes along the same lines as for theorem 2.3 of \cite{Thilliez2}, we shall skip some details. First, using the classical \L ojasiewicz inequality for analytic functions and the Cauchy formula, one obtains the existence of a neighborhood $ U $ of $ 0 $ in 
$ \mathbb{R}^n $ and constants $ C_1>0 $, $ \nu\geq 1 $, such that, for any multi-index $ J $ and any point $ x $ in $ U\setminus Z_g $, 
\begin{equation}\label{1overg}
\big\vert D^J\big(1/g\big)(x)\big\vert\leq C_1^{j+1}j!\dist(x,Z_g)^{-j-\nu}. 
\end{equation}
Now, let $ h $ be an element of $ \imax^\infty_{Y,M}$. One can then find a constant $ C_2>0 $ such that the inequality
\begin{equation}\label{hflat1} 
\vert D^Kh(x)\vert\leq
C_2^{k+1}k!M_k(C_2\dist(x,Y))^qM_q 
\end{equation}
holds for any multi-index $ K $, any integer $ q\geq 0 $ and any point $ x $ in $ U $: indeed, it suffices, maybe after shrinking $ U $, to apply the Taylor formula 
between $ x $ and a point $ \hat{x} $ of $ Y $ such that $ \vert x-\hat{x}\vert=\dist(x,Y) $. Taking the infimum with respect to $ q $ in \eqref{hflat1}, we obtain
\begin{equation}\label{hflat}
\vert D^Kh(x)\vert\leq C_2^{k+1}k!M_kh_M\big(C_2\dist(x,Y)\big).
\end{equation}
Now, using \eqref{hMsq} together with the definition of $ M^{(\theta)} $ recalled in subsection 
\ref{sequ} and the assumption \eqref{lojagen}, we have $ h_M\big(C_2\dist(x,Y)\big)\leq \big(h_M(C_3\dist(x,Y))\big)^2 $ and $ h_M\big(C_3\dist(x,Y)\big)\leq C_4h_{M^{(\theta)}}\big(C_4\theta(\dist(x,Y))\big)\leq C_4h_{M^{(\theta)}}\big(C_4\dist(x,Z_g)\big) $
for some suitable constants $ C_3 $ and $ C_4 $. This implies in particular, for any 
$ q\in \mathbb{N} $,
\begin{equation}\label{manip}
h_M\big(C_2\dist(x,Y)\big)\leq C_4^{q+1}M^{(\theta)}_q\big(\dist(x,Z_g)\big)^q h_M\big(C_3\dist(x,Y)\big).
\end{equation}
Now put $ L=J+K $, and note that we then have $ j+k=l $ and $ j!k!\leq l! $. Applying
\eqref{1overg}, \eqref{hflat} and \eqref{manip} with $ q=j+[\nu]+1 $, we obtain, for any $ x $ in $ U\setminus Z_g $, hence for any $ x $ in $ U\setminus Y $,
\begin{equation*}
\vert D^J\big(1/g\big)(x)D^Kh(x)\vert \leq C_5^{l+1}l!M_k M^{(\theta)}_{j+[\nu]+1} h_M\big(C_3\dist(x,Y)\big) 
\end{equation*}
with $ C_5=\max(C_1,C_4) $. Using remark \ref{grtMtheta} and properties \eqref{seq3} and \eqref{seq4} for the sequence $ M^{(\theta)} $,  
we also get $ M_k M^{(\theta)}_{j+[\nu]+1}\leq C_6^{l+1} M^{(\theta)}_l $ for some suitable $ C_6 $. 
For any multi-index $ L $, the Leibniz formula
\begin{equation*}
D^L(h/g)=\sum_{J+K=L}\frac{L!}{J!K!}D^J(1/g)D^Kh
\end{equation*}
and the preceding estimates then yield, for any point $ x $ in $ U\setminus Y $,
\begin{equation}\label{hoverg}
\vert D^L(h/g)(x)\vert \leq C_7^{l+1}l!M^{(\theta)}_l h_M\big(C_3\dist(x,Y)\big) 
\end{equation}
where $ C_7 $ is a constant. By \eqref{hoverg} and the Hestenes lemma, $ h/g $ extends to a $ C^\infty $ function $ \eta $ on $ U $; moreover $ \eta $ belongs to $ \imax^\infty_{Y,M^{(\theta)}} $. Since $ h=g\eta $ and $ g $ belongs to $ \mathcal{K} $, we therefore have shown
$ \imax^\infty_{Y,M}\subset \mathcal{K}\imax^\infty_{Y,M^{(\theta)}} $. Invoking
\eqref{dol}, we derive $ \langle\psi\rangle\imax^\infty_{Y,M}\subset \mathcal{K}\langle\psi\rangle\imax^\infty_{Y,M^{(\theta)}}\subset \langle\varphi\rangle\imax^\infty_{Y,M^{(\theta)}} $, which ends the proof.
\end{proof}

We end this section with an auxiliary lemma. 

\begin{lem}\label{null}
One has $ \imax^\infty_{Y,M}=\imax_M\imax^\infty_{Y,M} $.
\end{lem}

\begin{proof} Put $ g(x)=x_1^2+\cdots+x_n^2 $. It suffices to show that any element $ h $ of $ \imax^\infty_{Y,M} $ can be written $ h=g\eta $ with $ \eta\in  \imax^\infty_{Y,M} $. We can proceed as in the proof of theorem \ref{divis}: it is enough to remark that \eqref{1overg} is replaced by
$ \big\vert D^J\big(1/g\big)(x)\big\vert\leq C_7^{j+1}j!\vert x\vert^{-j-2} $ and
that the estimate \eqref{hflat1} for elements of $ \imax^\infty_{Y,M} $, applied with $ q=j+3 $, then implies
$ \big\vert D^J\big(1/g\big)(x)D^Kh(x)\big\vert\leq C_8^{j+k+1}k!M_{j+k}\vert x\vert $ for some positive constant $ C_8 $. The conclusion follows by the Leibniz formula as for theorem \ref{divis}. 
\end{proof}

\section{A Theorem of Infinite Determinacy}

\subsection{Analytic germs with non-isolated singularities}\label{fKf}
In what follows, we consider a real-analytic map germ $ \psi:(\mathbb{R}^n,0)\rightarrow(\mathbb{R}^p,0) $ and an element $ f $ of $ \int\langle\psi\rangle $. 
Since the Jacobian ideal $ \langle\nabla f\rangle $ is contained in $ \langle\psi\rangle $, we can use the machinery of subsection \ref{singuls} with $ \varphi=\nabla f $. The corresponding Fitting ideal $ \mathcal{K}$ will be denoted by $ \mathcal{K}_f $.

This framework has the following motivation. Let $ X $ be a real-analytic set germ at the origin and let $ \mathcal{I}_X $ be the ideal of elements of $ \mathcal{O}_n $ vanishing on $ X $. Let $ f $ be an element of $ \mathcal{O}_n $, with $ f(0)=0 $. The Bochnak-\L ojasiewicz inequality $ \vert f(x)\vert\leq C\vert x\vert\vert\nabla f(x)\vert $ (see \cite{Bochnak-Loj}) implies that $ f $ vanishes on its real critical set $ S_f $. Thus, $ S_f $ contains $ X $ if and only if $ f $ belongs to the primitive ideal $ \int\mathcal{I}_X $.

In order to study the determinacy properties of real-analytic germs whose critical locus contains the given set $ X $, we are therefore led to consider elements $ f$ of $ \int\langle\psi\rangle $ where $ \psi_1,\ldots,\psi_p $ is a system of generators of $ \mathcal{I}_X $.

\subsection{The main result}
Let $ \psi $, $ f $ and $ \mathcal{K}_f $ be as in subsection \ref{fKf} above. We then have the following statement.

\begin{thm}\label{main}
Let $ Y $ be a germ of closed set at the origin of $ \mathbb{R}^n $ and let $ \theta $ be an admissible function. Assume that the ideal $ \mathcal{K}_f $ contains a germ $ g $ whose complex zero set $ Z_g $ satisfies the generalized \L ojasiewicz inequality
\begin{equation*}
\dist(x,Z_g)\geq \theta\big(\dist(x,Y)\big)\text{ for any real point }x.
\end{equation*}
Then, for any non-quasianalytic tame sequence $ M $, we have 
\begin{equation*}
f+\int\langle\psi\rangle\imax^\infty_{Y,M}\subset f\circ\mathcal{R}_{M^{(\theta)}}.
\end{equation*}
\end{thm}
\begin{proof} We use the classical homotopy method of Mather. Let $ u $ be an element of $ \int\langle\psi\rangle\imax^\infty_{Y,M} $. 
For any $ x $ in a sufficiently small neighborhood $ V $ of $ 0 $ and any $ t $ in $ [0,1] $, put $ \widetilde{f}(x,t)=f(x)+tu(x) $. The key of the proof now lies in the following claim. 
\begin{claim} 
For any $ a $ in $ [0,1] $, one can find a
time-dependent vector field $ \Xi_a $ of $ \mathbb{R}^n $, whose components belong to $ \imax^\infty_{Y,M^{(\theta)}}(a) $ (see subsection \ref{classes} for notation), 
and such that the equality
\begin{equation}\label{cruc}
\frac{\partial\widetilde{f}}{\partial t}(x,t)=\left\langle\nabla_x\widetilde{f}(x,t)\,,\,
\Xi_a(x,t)\right\rangle 
\end{equation}
holds for any $ x $ in a neighborhood of $ 0 $ in $ \mathbb{R}^n $ and any $ t $ in a neighborhood of $ a $ in $ [0,1] $ (the brackets $ \langle\cdot,\cdot\rangle $ denote the usual euclidean scalar product in $ \mathbb{R}^n $).
\end{claim}
Admitting temporarily the claim, the rest of the proof follows a standard scheme. For the reader's convenience, we briefly recall the argument. Note first that ordinary differential equations can be solved in ultradifferentiable classes, see e.g. \cite{Dynkin}. Thus, integrating the equation
\begin{equation}\label{EDO1}
\frac{\partial\phi}{\partial t}(\xi,t)=-\Xi_a\big(\phi(\xi,t),t\big)
\end{equation}
with the initial condition
\begin{equation}\label{EDO2}
\phi(\xi,a)=\xi
\end{equation}
yields a solution $ \phi_a $ defined in a neighborhood of $ (0,a) $ in $ \mathbb{R}^n\times [0,1] $, whose components belong to $ \mathcal{C}_{n,M^{(\theta)}}(a) $, and such that
for any $ t $ sufficiently close to $ a $, the map $ \xi\mapsto\phi_a(\xi,t) $ is a germ of diffeomorphism at the origin of $ \mathbb{R}^n $.
Using \eqref{cruc}, \eqref{EDO1} and \eqref{EDO2}, it is then easy to obtain
\begin{equation*}
\widetilde{f}\big(\phi_a(\xi,t),t\big)=\widetilde{f}(\xi,a)
\end{equation*}
for any $ (\xi,t) $ in a neighborhood of $ (0,a) $. By an immediate compactness argument,
there exists a finite family of points $ 0=a_0<\cdots<a_m=1 $ such that 
$ \widetilde{f}\big(\phi_{a_j}(\xi,a_{j+1}),a_{j+1}\big)=\widetilde{f}(\xi,a_j) $ for $ j=0,\ldots, m-1 $. Consider the germ of diffeomorphism $ \Phi=\big(\phi_{a_0}(\,\cdot\,,a_1)\big)^{-1}\circ\cdots\circ\big(\phi_{a_{m-1}}(\,\cdot\,,a_m)\big)^{-1} $. We have $ f+u=f\circ\Phi $ by construction. By standard results on stability of ultradifferentiable classes under composition and inversion of mappings \cite{Dynkin}, the components of $ \Phi $ belong to $ \mathcal{C}_{n,M^{(\theta)}} $,  
hence the theorem.

Thus, it remains to proving the claim above. Note first that the germ $ u $ belongs to
$ \langle\nabla f\rangle\imax^\infty_{Y,M^{(\theta)}} $ by virtue of theorem \ref{divis}.
Taking remark \ref{injec} into account, we get, for any $ a $ in $ [0,1] $, 
\begin{equation}\label{EDO3}
\frac{\partial\tilde{f}}{\partial t}\in 
\langle\nabla f\rangle\imax^\infty_{Y,M^{(\theta)}}(a).
\end{equation}
We also have  
\begin{equation*}
\langle\nabla f\rangle\imax^\infty_{Y,M^{(\theta)}}(a)\subset 
\langle\nabla_x\widetilde{f}\rangle\imax^\infty_{Y,M^{(\theta)}}(a)+\langle\nabla u\rangle\imax^\infty_{Y,M^{(\theta)}}(a).
\end{equation*}
Moreover, the derivatives of $ u $ also belong to $ \langle\psi\rangle\imax^\infty_{Y,M} $, 
hence to $ \langle\nabla f\rangle\imax^\infty_{Y,M^{(\theta)}} $ by theorem \ref{divis}.
Using lemma \ref{null} (with $ M^{(\theta)} $ instead of $ M $), we have $ \imax^\infty_{Y,M^{(\theta)}}\subset \imax_{M^{(\theta)}}\imax^\infty_{Y,M^{(\theta)}}\subset
\imax_{M^{(\theta)}}(a)\imax^\infty_{Y,M^{(\theta)}}(a)  $.
Thus, we obtain 
\begin{equation*}
\langle\nabla f\rangle\imax^\infty_{Y,M^{(\theta)}}(a)\subset 
\langle\nabla_x\tilde{f}\rangle\imax^\infty_{Y,M^{(\theta)}}(a)+\imax_{M^{(\theta)}}(a)\langle\nabla f\rangle
\imax^\infty_{Y,M^{(\theta)}}(a).
\end{equation*}
By Nakayama's lemma, this implies
$\langle\nabla f\rangle\imax^\infty_{Y,M^{(\theta)}}(a)\subset 
\langle\nabla_x\tilde{f}\rangle\imax^\infty_{Y,M^{(\theta)}}(a) $, 
which, together with \eqref{EDO3}, yields the claim. The proof is complete.
\end{proof}

\begin{rem}\label{presX}
An inspection of the preceding proof shows that it provides a diffeomorphism $ \Phi $ that coincides with the identity at infinite order on $ Y $. Beside this, the equality 
$ f+u=f\circ\Phi $ with both $ f $ and $ u $ singular on $ X=\psi^{-1}(\{0\}) $ implies
$ \Phi(X)\subset S_f $, where we recall that $ S_f $ denotes the critical set of $ f $. Thus, if $ X = S_f $ (the assumption $ f\in\int\langle\psi\rangle $ only implies the inclusion), we derive that $ \Phi $ preserves $ X $. By a slight modification of the proof of theorem \ref{main}, one can show that it is also the case when $ \int\langle\psi\rangle=\langle\psi\rangle^2 $. 
\end{rem}

\begin{rem}\label{typic} 
Being given an element $ g$ of $ \mathcal{K}_f $, theorem \ref{main} can always be applied with $ Y= Z_g\cap\mathbb{R}^n $, or more generally with $ Y $ chosen as a subanalytic set containing $ Z_g\cap\mathbb{R}^n $. Indeed, in a neighborhood $ 0 $ in $ \mathbb{R}^n,
$ both functions $ x\mapsto\dist(x,Z_g) $ and $ x\mapsto\dist(x,Y) $ are then
subanalytic and the second one 
vanishes at any point where the first one vanishes. Hence inequality \eqref{lojagen} holds (and is sharp) with $ \theta(t)=ct^s, $ where $ c $ is a
suitable positive constant and the \L ojasiewicz exponent $ s $ is a rational number (see e.g. \cite{Bochnak-Risler}).
\end{rem}

\begin{rem}\label{smoothcase}
The proof of theorem \ref{main} can be easily adapted to the $ C^\infty $ case (it is then simpler) in order to obtain the implication \eqref{jacX0}$\Longrightarrow$\eqref{deterX0} mentioned in the introduction. 
\end{rem}

\subsection{The particular case of isolated singularities}\label{isolsi}
From theorem \ref{main}, it is possible to recover theorem 4.1 of \cite{Thilliez1} as follows. Put $ p=n $ and $ \psi_j(x)=x_j $ for $ j=1,\ldots,n $, so that $ \int\langle\psi\rangle=\langle x\rangle^2 $. The relations $ k^1,\ldots,k^q $ of subsection \ref{singuls} are the trivial ones: $ x_ke_l-x_le_k $ for $ 1\leq k<l\leq n $, where the $ e_i $ are the elements of the canonical basis of $ \mathcal{O}_n^n $. It is then not difficult to check that $ \mathcal{K}_f $ contains all the germs $ x_k^{n-2}\partial f/\partial x_j $, hence it contains $ \vert x\vert^{2n}\partial f/\partial x_j $ for $ j=1,\ldots,n $. Now, as in theorem 4.1 of \cite{Thilliez1}, assume that $ \langle\nabla f\rangle $ contains an element $ \gamma $ whose complex zero set $ Z_\gamma $ satisfies
$ d(x,Z_\gamma)\geq c\vert x\vert^s $ for real points $ x $, with $ c>0 $ and $ s\geq 1 $. Then the germ $ g=\vert x\vert^{2n}\gamma $ belongs to $ \mathcal{K}_f $ and we have  $ Z_g = Z_\gamma\cup\Gamma $, with $ \Gamma=\{x\in\mathbb{C}^n: x_1^2+\cdots+x_n^2=0\} $. Since $ \Gamma $ and $ \mathbb{R}^n $ satisfy the regular separation property with 
\L ojasiewicz exponent $ 1 $, we derive $ d(x,Z_g)\geq c\vert x\vert^s $ for real points $ x $, maybe after changing $ c $ (but not $ s$).  
Thus, the assumption of theorem \ref{main} holds with $ Y=\{0\} $ and $ \theta(t)=ct^s $. 
Since, in this case, we also have $ \int\langle\psi\rangle\imax^\infty_{Y,M}=\imax^\infty_M $ and 
$ M^{(\theta)}=M^s $ (see example \ref{expow}), we 
obtain finally $ f+\imax^\infty_M\subset f\circ\mathcal{R}_{M^s} $, which is precisely the conclusion of theorem 4.1 of \cite{Thilliez1}.

\subsection{Examples}\label{examp}

\begin{exam}\label{exisol}
The following example, mentioned in the introduction, illustrates the influence of flatness on $ Y $, even in the case of isolated singularities. Put $ n=2 $ and $ f(x)=(x_1^2+x_2^4)^2 $. Being given positive real numbers $ \alpha $ and $ \mu $, we study the determinacy of $ f $  
with respect to Gevrey perturbations (that is, with $ M_j=j!^\alpha $), for which we assume flatness on $ Y=\{ x\in\mathbb{R}^2 :\vert x_1\vert =\vert x_2\vert^\mu\} $. Proceeding as in subsection \ref{isolsi}, we see that the ideal $ \mathcal{K}_f $ is generated by $ x_1(x_1^2+x_2^4) $ and $ x_2^3(x_1^2+x_2^4) $; in particular, $ f $ belongs to $ \mathcal{K}_f $. The zero set $ Z_f $ is the union of the two smooth curves $ z_1=\pm i z_2^2 $. We derive $ \dist(x,Z_f)\approx\min\big(\vert x_1-ix_2^2\vert,\vert x_1+ix_2^2\vert\big)\approx \vert x_1\vert+x_2^2 $ for $ x $ real.  
For $ \mu\geq 1 $, we have $ \dist(x,Y)\approx \min\big(\vert x_1-\vert x_2\vert^\mu\vert,\vert x_1+\vert x_2\vert^\mu\vert\big) $, hence $ \dist(x,Y)\lesssim \vert x_1\vert+\vert x_2\vert^\mu $. This yields $ \theta(\dist(x,Y))\leq \dist(x,Z_f) $ with $ \theta(t)\approx t^{2/\mu} $ in the case $ 1\leq\mu\leq 2 $, and
$ \theta(t)\approx t $ in the case $ \mu>2 $. For $ 0<\mu<1 $, similar arguments give $ \dist(x,Y)\lesssim \vert x_1\vert^{1/\mu}+\vert x_2\vert $, hence $ \theta(\dist(x,Y))\leq \dist(x,Z_f) $ with $ \theta(t)\approx t^2 $.
Since we have $ \int\langle\psi\rangle=\langle x\rangle^2 $, lemma \ref{null} yields $ \int\langle\psi\rangle\imax^\infty_{Y,M}=\imax^\infty_{Y,M} $. Applying theorem \ref{main} and taking into account example \ref{expow}, we finally obtain, with the notation of the introduction,
$ f+\imax^\infty_{Y,\alpha}\subset f\circ\mathcal{R}_\beta $ with 
\begin{equation}\label{casebeta}
\beta=\left\{
\begin{array}{ll}
2\alpha & \text{ if }0<\mu\leq 1,\\
2\alpha/\mu & \text{ if }1<\mu\leq 2,\\
\alpha & \text{ if }\mu>2.
\end{array}
\right.
\end{equation}
Similar computations show that we can take $ \beta =\alpha $ when $ Y $ is the 
$ x_2 $-axis and $ \beta=2\alpha $ when $ Y $ is the $ x_1 $-axis.
We know from \cite{Thilliez1} that the result is sharp. 
\end{exam}

\begin{exam} Here is the other example mentioned in the introduction. Put $ n=3 $ and
$ f(x)=(x_1^2+x_2^2+x_3^4)^2(x_1^2+x_2^2) $. The critical locus of $ f $ is the $ x_3 $-axis, and we can take $ \psi(x)=(x_1,x_2) $. We then have $ \int\langle\psi\rangle=\langle x_1^2,x_1x_2,x_2^2\rangle $. The ideal $ \mathcal{K}_f $ contains the germ $ g(x)=(x_1^2+x_2^2+x_3^4)^2(3x_1^2+3x_2^2+x_3^4)^2 $. Put $ 
Y=\{x\in\mathbb{R}^3 : \vert x_1\vert=\vert x_2\vert=\vert x_3\vert^\mu\} $. Using the same kind of computations as in example \ref{exisol}, it is not difficult to obtain $ \theta(\dist(x,Y))\leq \dist(x,Z_g) $, where $ \theta $ depends on $ \mu $ in the same way (remark that, since the complex zero set of the ideal $ \mathcal{K}_f $ is precisely given by $ x_1^2+x_2^2+x_3^4=0 $, one cannot expect to find any $ g $ in $ \mathcal{K}_f $ providing better estimates). In the case of Gevrey perturbations, theorem \ref{main} yields $ f+\langle x_1^2,x_1x_2,x_2^2\rangle\imax^\infty_{Y,\alpha}\subset f\circ\mathcal{R}_\beta $,
where $\beta $ is given by \eqref{casebeta}. As mentioned in the introduction, we are not limited to Gevrey sequences: for instance, if $ Y=\{x\in\mathbb{R}^3 : \vert x_1\vert=\vert x_2\vert=\vert x_3\vert^\mu\Log(1+\frac{1}{\vert x_3\vert})^{\nu} \} $ with $ 1<\mu<2 $ and $ \nu>0 $, we obtain 
$ f+\langle x_1^2,x_1x_2,x_2^2\rangle\imax^\infty_{Y,\alpha}\subset
f\circ\mathcal{R}_{M^+} $ with $ M^+_j=j!^{2\alpha/\mu}(\Log j)^{2\nu j/\mu} $.
This corresponds to the application of theorem \ref{main} with 
$ \theta(t)\approx t^{2/\mu}\Log(1+\frac{1}{t})^{-\nu} $, taking into account example \ref{expowlog}.
\end{exam}

\subsection{A geometric question}
Denote by $ \Delta_f $ the complex zero set of $ \mathcal{K}_f $. It would be interesting to obtain a sufficient condition of determinacy involving only a \L ojasiewicz separation property between $ \Delta_f $ and $ Y $, instead of the choice of a particular germ $ g $ in $ \mathcal{K}_f $. Thus, we are led to the following question.

\begin{prob}
Let $ \Delta $ be a germ of complex analytic set at the origin of $ \mathbb{C}^n $ and assume that $ \Delta\cap\mathbb{R}^n=\{0\} $. Can one find a purely $ (n-1) $-dimensional complex analytic set $ Z $ containing $ \Delta $, such that $ Z\cap\mathbb{R}^n=\{0\} $ and
$ \dist(x,Z)\approx\dist(x,\Delta) $ for any $ x $ in a neighborhood of $ 0 $ in
$ \mathbb{R}^n $?
\end{prob}

A positive answer could be applied as follows: take $ \Delta=\Delta_f $ and assume that such a $ Z $ exists. Then there is a holomorphic function germ $ h $ such that $ Z=h^{-1}(\{0\}) $. Since the function $ h $ vanishes on $ \Delta_f $, the Nullstellensatz yields an integer $ m\geq 1 $ such that $ h^m $ belongs to the complexification of $ \mathcal{K}_f $. For $z\in\mathbb{C}^n $, put $ g(z)=h^m(z)\overline{h^m(\overline{z})} $. It is then easy to see that $ g $ belongs to $ \mathcal{K}_f $ and that $ Z_g= Z\cup S(Z) $, where $ S $ denotes the conjugation map $ z\mapsto \overline{z} $. Hence, for real points $ x $, we have easily
$ \dist(x,Z_g)=\dist(x,Z)\approx\dist(x,\Delta_f) $. The \L ojasiewicz-type condition
$ \dist(x,\Delta_f)\geq \theta(\dist(x,Y)) $ would then suffice to apply theorem \ref{main}.


\begin{thebibliography}{00}

\bibitem{Bochnak-Loj} \textsc{J. Bochnak \& S. \L ojasiewicz}, \textit{A converse of the Kuiper-Kuo theorem}, Proceedings of Liverpool Singularities, Springer Lecture Notes in Math. {\bf 192} (1971), 254--261. 

\bibitem{Bochnak-Risler} \textsc{J. Bochnak \& J.-J. Risler}, \textit{Sur les exposants de \L ojasiewicz}, Comment. Math. Helvetici \textbf{50} (1975), 493--507.

\bibitem{Dynkin} \textsc{E.M. Dynkin}, \textit{Pseudoanalytic extension of smooth functions. The uniform scale}, Amer. Math. Soc. Transl. \textbf{115} (1980), 33--58.

\bibitem{Eisenbud} \textsc{D. Eisenbud}, \textit{Commutative Algebra with a view toward Algebraic Geometry}, Graduate Texts in Math. 150, Springer Verlag (1995).

\bibitem{Grandjean2} \textsc{V. Grandjean}, \textit{Infinite relative determinacy of smooth function germs with transverse isolated singularities and relative \L ojasiewicz conditions}, J. London Math. Soc. \textbf{69} (2004), 518--530.

\bibitem{Kushner-Terra Leme} \textsc{L. Kushner \& B. Terra Leme}, \textit{Finite relative determination and relative stability}, Pacific J. Math. \textbf{192} (2000), 315--328.

\bibitem{Nguyen} \textsc{T.C. Nguyen, H.D. Nguyen, S.M. Nguyen \& H.V. Ha}, \textit{Sur les germes de fonctions infiniment d\'etermin\'es}, C. R. Acad. Sci. Paris  \textbf{285} (1977), 1045--1048.

\bibitem{Sun-Wilson} \textsc{B. Sun \& L.C. Wilson}, \textit{Determinacy of smooth germs with real isolated line singularities}, Proc. Amer. Math. Soc. \textbf{129} (2001), 2789--2797.

\bibitem{Thilliez0} \textsc{V. Thilliez}, \textit{Sur les fonctions compos\'ees ultradiff\'erentiables}, J. Math. Pures et Appl. \textbf{76} (1997), 499--524.

\bibitem{Thilliez1} \textsc{V. Thilliez}, \textit{Germes de d\'etermination infinie dans des classes lisses}, Manuscripta Math. \textbf{99} (1999), 203--222.

\bibitem{Thilliez2} \textsc{V. Thilliez}, \textit{On closed ideals in smooth classes}, Math. Nachr. \textbf{227} (2001), 143--157.

\bibitem{Thilliez3} \textsc{V. Thilliez}, \textit{Infinite determinacy on a closed set for smooth germs with non-isolated singularities}, Proc. Amer. Math. Soc. (to appear).

\bibitem{Wall} \textsc{C.T.C. Wall}, \textit{Finite determinacy of smooth map-germs}, Bull. London Math. Soc. \textbf{13} (1981), 481--539.

\bibitem{Wilson} \textsc{L.C. Wilson}, \textit{Infinitely determined mapgerms}, Can. J. Math. \textbf{33} (1981), 671--684.

\end{thebibliography}
\end{document}